\documentclass[preprint,11pt]{elsarticle}
\usepackage[T1]{fontenc}
\usepackage{fullpage}
\usepackage{hyperref}
\usepackage{amssymb}
\usepackage{amsthm,bm,mathtools}
\usepackage{lineno}
\usepackage{mathrsfs}
\usepackage{amsmath}
\usepackage{color}
\usepackage{xcolor}
\hypersetup{colorlinks=true,linkcolor=blue!55!black,urlcolor=blue!55!black}
\numberwithin{equation}{section}
\newtheorem{theorem}{Theorem}[section]
\newtheorem{lemma}[theorem]{Lemma}
\newtheorem{proposition}[theorem]{Proposition}
\newtheorem{corollary}[theorem]{Corollary}
\newtheorem{conjecture}[theorem]{Conjecture}
\theoremstyle{definition}
\newtheorem{definition}[theorem]{Definition}
\newtheorem{example}[theorem]{Example}

\newcommand{\rk}{\operatorname{rk}}
\newcommand{\cl}{\operatorname{cl}}
\newcommand{\Des}{\operatorname{Des}}
\newcommand{\cH}{\mathcal H}
\newcommand{\Hd}{H^{*}}
\newcommand{\Fd}{F^{*}}

\newcommand{\U}{\mathsf U}
\newcommand{\Gt}{\widetilde\Gamma}
\newcommand{\des}{\mathrm{des}}
\newcommand{\zmin}{\widehat 0}
\newcommand{\zmax}{\widehat 1}

\allowdisplaybreaks[4]
\begin{document}
\begin{frontmatter}
\title{Real-rootedness of dual Chow polynomials of paving matroids}
\author{Yanxin Liu}
\ead{liu-yanxin@outlook.com}

\author{Lijie Wang}
\ead{lijiewang26@hotmail.com}

\author{Qiqi Xiao\corref{cor1}}
\ead{xiaoqqcs@hotmail.com}
\cortext[cor1]{Corresponding author}

\address{School of Mathematical Sciences, Dalian University of Technology, Dalian 116024, P.R. China}
\begin{abstract}
We prove that the dual Chow polynomial and the dual right augmented Chow polynomial of every finite paving matroid have only negative real zeros. Moreover, the dual Chow polynomial interlaces the dual right augmented Chow polynomial. This confirms the conjectures of Caiolo, Ferroni, and Hoster for paving matroids.
\end{abstract}
\begin{keyword}
Paving matroids \sep
Dual Chow polynomials \sep
Dual right augmented Chow polynomials \sep
Real-rootedness 
\MSC[2020] 05B35 \sep 26C10 \sep 06A07
\end{keyword}
\end{frontmatter}

\section{Introduction}
The Chow ring of a loopless matroid is a graded algebra
constructed from its lattice of flats by Feichtner and
Yuzvinsky~\cite{FY04}.
Its Hilbert--Poincar\'e series is called the Chow polynomial of the matroid.
For a matroid of rank $r$, Poincar\'e duality implies that
the Chow polynomial is symmetric with center of symmetry
$(r-1)/2$, while the hard Lefschetz theorem further implies
that its coefficients are unimodal~\cite{AHK18}.
Recall that a matroid of rank $r$ is paving if every circuit has at least $r$ elements; see~\cite{Oxl11}.
Br\"and\'en and Vecchi~\cite{BV26} proved that the Chow and
augmented Chow polynomials of every paving matroid are
real-rooted.
For general matroids, however, B\'erczi, Coron, Eur, and Li~\cite{BCEL26} constructed counterexamples to log-concavity, and hence to real-rootedness, for both polynomials.

Ferroni, Matherne, and Vecchi~\cite{FMV24} developed a general theory of Chow functions for partially ordered sets, extending the study of Chow polynomials beyond the setting of Chow rings.
Building on this framework, Caiolo, Ferroni, and Hoster~\cite{CFH26} introduced the dual Chow polynomial $\Hd_M(x)$ and the dual right augmented Chow polynomial $\Fd_M(x)$ of a matroid $M$ by considering the sign-twisted reverse of the characteristic kernel of its lattice of flats.
Unlike the ordinary Chow polynomial, whose constant term is one,
the dual Chow polynomial has constant term equal to the absolute
value of the M\"obius function of the lattice of flats of $M$,
evaluated at its minimum and maximum elements.
The dual Chow polynomials are also related to the inverse Kazhdan--Lusztig polynomial~\cite{GX21} and the inverse $Z$-polynomial~\cite{FMSV24,GRX26} of a matroid (see~\cite[Section~7.3]{CFH26}).

A geometric motivation for studying dual Chow polynomials
comes from their connection with Orlik--Solomon algebras.
Based on this connection, Caiolo, Ferroni, and Hoster~\cite[Section~7.1]{CFH26} raised the question of whether the dual Chow polynomial can be realized as the Hilbert--Poincar\'e series of a natural graded module over the Chow ring and whether such a module admits a Lefschetz package.
Although these questions remain open, several combinatorial properties of dual Chow polynomials have been established.
In particular, Caiolo, Ferroni, and Hoster~\cite{CFH26} obtained chain formulas, deletion formulas, and explicit expansions in terms of the flag $h$-vector.
The latter imply that both $\Hd_M(x)$ and $\Fd_M(x)$ are $\gamma$-positive, with centers of symmetry $(r-1)/2$ and $r/2$, respectively, for a matroid of rank $r$.

Despite these positivity properties, real-rootedness does not hold for dual Chow polynomials of general Cohen--Macaulay posets, as shown by a counterexample in \cite[Remark~7.4]{CFH26}.
For matroids, however, Caiolo, Ferroni, and Hoster~\cite[Corollary~6.11]{CFH26} proved that both $\Hd_M(x)$ and $\Fd_M(x)$ are real-rooted when $M$ is a uniform matroid and proposed the following conjecture.

\begin{conjecture}[{\cite[Conjecture~7.3]{CFH26}}]
For every matroid $M$, the dual Chow polynomial $\Hd_M(x)$
has only nonpositive real zeros.
\end{conjecture}

They further conjectured that $\Fd_M(x)$ has only nonpositive real zeros and that $\Hd_M(x)$ interlaces $\Fd_M(x)$.
In this paper, we establish these conjectured properties for paving matroids.
Our main result is the following.

\begin{theorem}\label{thm:main}
Let $M$ be a finite paving matroid.
Then both $\Hd_M(x)$ and $\Fd_M(x)$ have only negative real zeros, and $\Hd_M(x)$ interlaces $\Fd_M(x)$.
\end{theorem}

The paper is organized as follows.
In Section~\ref{sec:pre}, we recall the necessary definitions and properties of dual Chow polynomials, together with several results on interlacing polynomials.
In Section~\ref{sec:decomposition}, we establish recurrence relations for the dual Chow and dual right augmented Chow polynomials of paving matroids.
In Section~\ref{sec:eulerian}, we use these identities and the interlacing properties of the refined Eulerian polynomials to prove Theorem~\ref{thm:main}.

\section{Preliminaries}\label{sec:pre}
In this section, we introduce the notation and terminology concerning matroids and posets, recall the basic properties of dual Chow polynomials, and collect several results on interlacing polynomials that will be used in our proofs.
\subsection{Matroids and posets}
We begin with some basic notions concerning matroids.
For further background on matroid theory, we refer the
reader to Oxley~\cite{Oxl11}.

\begin{definition}
A \emph{matroid} is a pair $M=(E,\mathcal I)$ consisting of a finite set $E$, called the \emph{ground set}, and a collection $\mathcal I$ of subsets of $E$, called \emph{independent sets}, satisfying the following properties:
\begin{enumerate}
\item[(i)] $\varnothing\in\mathcal I$;
\item[(ii)] If $I\in\mathcal I$ and $J\subseteq I$,
then $J\in\mathcal I$;
\item[(iii)] If $I,J\in\mathcal I$ and $|I|<|J|$,
then there exists $e\in J\setminus I$ such that
$I\cup\{e\}\in\mathcal I$.
\end{enumerate}
\end{definition}

A subset of $E$ that is not independent is called \emph{dependent}.
A maximal independent set is called a \emph{basis}, and a minimal dependent set is called a \emph{circuit}.
The \emph{rank function} of $M$ is defined by
\[
\rk_M(A)=\max\{|I|:I\subseteq A,\ I\in\mathcal I\},
\quad A\subseteq E.
\]
The \emph{rank} of $M$ is $\rk(M)=\rk_M(E)$.
When the underlying matroid is clear from the context, we simply write $\rk(A)$ for $\rk_M(A)$.
An element $e\in E$ is called a \emph{loop} if $\{e\}$ is a circuit.
A matroid is \emph{loopless} if it has no loops.
A subset $F\subseteq E$ is called a \emph{flat} if $\rk(F\cup\{e\})>\rk(F)$ for every $e\in E\setminus F$.
The \emph{closure} $\cl(A)$ of a subset $A\subseteq E$ is the smallest flat containing $A$.
The flats of $M$, ordered by inclusion, form a lattice, called the \emph{lattice of flats} of $M$ and denoted by $\mathcal L(M)$.
Its minimum and maximum elements are $\cl(\varnothing)$ and $E$, respectively.
In particular, if $M$ is loopless, then the minimum element of $\mathcal L(M)$ is $\varnothing$.
A flat distinct from $E$ is called a \emph{proper flat}, and a maximal proper flat is called a \emph{hyperplane}.
Equivalently, a hyperplane of a rank-$r$ matroid is a flat of rank $r-1$.

For integers \(0\leq r\leq n\), the \emph{uniform matroid} \(\mathrm{U}_{r,n}\) on an n-element ground set has independent sets exactly the subsets of size at most $r$.
Recall that a matroid of rank $r$ is paving if every circuit has at least $r$ elements, or equivalently, if every subset of cardinality at most $r-1$ is independent.
In particular, every uniform matroid is paving.

Since dual Chow polynomials are associated with lattices of flats, we next recall some basic notions for posets.
For background on posets and incidence algebras, see Stanley~\cite [Chapter~3]{Sta12}.
A finite poset $P$ is called \emph{bounded} if it has a minimum
element $\zmin$ and a maximum element $\zmax$.
For $s\leq t$ in $P$, the \emph{interval} $[s,t]$ is defined by
\[
[s,t]=\{w\in P:s\leq w\leq t\}.
\]
A finite bounded poset is called \emph{graded} if all its maximal chains have the same length, where the length of a chain is the number of strict inequalities in the chain.
For a graded poset $P$, the \emph{rank} $\rk(t)$ of an element $t\in P$ is the length of a maximal chain in $[\zmin,t]$, and the rank of $P$ is $\rk(\zmax)$.
In particular, the lattice of flats of a matroid is graded, and the rank of a flat in $\mathcal L(M)$ coincides with its rank in $M$.

The \emph{M\"obius function} $\mu_P$ of a finite poset $P$ is determined recursively by
\[
\mu_P(s,s)=1,\qquad
\mu_P(s,t)=-\sum_{s\leq w<t}\mu_P(s,w)
\quad\text{if }s<t.
\]
For $P=\mathcal L(M)$, we write $\mu_M$ for $\mu_P$.
A \emph{Boolean lattice} is a poset isomorphic to the lattice of subsets of a finite set ordered by inclusion.
If $P$ is a Boolean lattice, then $\mu_P(s,t)=(-1)^{\rk(t)-\rk(s)}$ for every interval $[s,t]$ of $P$.

We next recall the notions of flag $f$-vectors and flag $h$-vectors (see~\cite[Section~3.13]{Sta12}).
Let $P$ be a finite bounded graded poset of rank $r\geq 1$.
Throughout, we write $[n]=\{1,\ldots,n\}$ for a positive integer $n$ and set $[0]=\varnothing$.
Given $S\subseteq[r-1]$, the \emph{rank-selected subposet} of $P$ on $S$ consists of the elements whose ranks belong to $S$, with the order inherited from $P$.
Writing $S=\{s_1<\cdots<s_k\}$, we define
\begin{align*}
\alpha_P(S)
&=\left|\left\{
t_1<\cdots<t_k:\ t_i\in P,\ \rk(t_i)=s_i,\ 
1\leq i\leq k
\right\}\right|,\\
\beta_P(S)
&=\sum_{T\subseteq S}(-1)^{|S\setminus T|}\alpha_P(T),
\end{align*}
with $\alpha_P(\varnothing)=1$.
The \emph{flag $f$-vector} and the \emph{flag $h$-vector} of $P$ are given by $(\alpha_P(S))_{S\subseteq[r-1]}$ and $(\beta_P(S))_{S\subseteq[r-1]}$, respectively.
For a matroid $M$, we write $\beta_M(S)$ for $\beta_{\mathcal L(M)}(S)$.

\subsection{Dual Chow polynomials and $\gamma$-expansions}
Let $P$ be a finite bounded graded poset.
The \emph{characteristic kernel} $\chi$ of $P$ is given by
\[
\chi_{st}(x)=
\sum_{s\leq w\leq t}
\mu_P(s,w)x^{\rk(t)-\rk(w)}.
\]
By considering the sign-twisted reverse of the characteristic kernel, Caiolo, Ferroni, and Hoster~\cite{CFH26} introduced the \emph{dual Chow polynomial} $\Hd_P(x)$ and the \emph{dual right augmented Chow polynomial} $\Fd_P(x)$.
We refer to~\cite[Section~3]{CFH26} for their original definitions and recall below their expressions in terms of the flag $h$-vector.

For a matroid $M$ with lattice of flats $\mathcal L(M)$,
we write
\[
\Hd_M(x)=\Hd_{\mathcal L(M)}(x)
\quad\text{and}\quad
\Fd_M(x)=\Fd_{\mathcal L(M)}(x).
\]
For a one-element poset $P$, we have
$\Hd_P(x)=\Fd_P(x)=1$.

Before presenting the flag $h$-vector expressions, we recall the notions of symmetry and $\gamma$-positivity.
For a survey on $\gamma$-positivity in combinatorics and geometry, see Athanasiadis~\cite{Ath18}.
Let $d$ be a nonnegative integer.
A polynomial $f(x)=\sum_{i}a_ix^i$ is said to be \emph{symmetric with center of symmetry $d/2$} if $a_i=a_{d-i}$ for every $i\in\mathbb Z$, where $a_i=0$ for negative values of $i$.
Equivalently, $x^df(x^{-1})=f(x)$.
If the center of symmetry is $d/2$ for $d=\deg f(x)$, the polynomial $f(x)$ will be called \emph{palindromic}.
A palindromic polynomial with center of symmetry $d/2$ admits a unique expansion of the form
\[
f(x)=
\sum_{i=0}^{\lfloor d/2\rfloor}
\gamma_i x^i(1+x)^{d-2i},
\]
called its \emph{$\gamma$-expansion}.
The polynomial $f(x)$ is said to be \emph{$\gamma$-positive} if $\gamma_i\geq 0$ for every $i$.
The associated \emph{$\gamma$-polynomial} is defined by
\[
\gamma(f;z)=
\sum_{i=0}^{\lfloor d/2\rfloor}\gamma_i z^i.
\]
Equivalently,
\begin{equation*}
f(x)=(1+x)^d
\gamma\left(f;\frac{x}{(1+x)^2}\right).
\end{equation*}
In particular, $\gamma(f;0)=f(0)$.
For the zero polynomial, we set $\gamma(0;z)=0$.
Note that $\gamma$-positivity implies symmetry and unimodality.

A subset of $\mathbb Z$ is called \emph{stable} if it
contains no two consecutive integers.
For $S\subseteq[r-1]$, we write $S^c=[r-1]\setminus S$.
Caiolo, Ferroni, and Hoster~\cite{CFH26} obtained explicit $\gamma$-expansions for the dual Chow and dual right augmented Chow polynomials.
More precisely, for a finite bounded graded poset $P$ of rank $r\geq 1$, we have
\begin{equation}\label{eq:intro-flag-H}
\Hd_P(x)=
\sum_{\substack{
S\subseteq[r-1]\\
r-1\notin S\\
S\text{ stable}
}}
\beta_P(S^c)x^{|S|}(1+x)^{r-1-2|S|}
\end{equation}
and
\begin{equation}\label{eq:intro-flag-F}
\Fd_P(x)=
\sum_{\substack{
S\subseteq[r-1]\\
S\text{ stable}
}}
\beta_P(S^c)x^{|S|}(1+x)^{r-2|S|}.
\end{equation}
These expansions show that $\Hd_P(x)$ and $\Fd_P(x)$ are symmetric about $(r-1)/2$ and $r/2$, respectively.
Moreover, both polynomials are $\gamma$-positive whenever the flag $h$-vector of $P$ is nonnegative.
Since the lattice of flats of a matroid is Cohen--Macaulay and hence has a nonnegative flag $h$-vector (see~\cite{Bjo92}), both $\Hd_M(x)$ and $\Fd_M(x)$ are $\gamma$-positive for every matroid $M$.

\subsection{Interlacing polynomials}
We next recall several properties of interlacing
polynomials that will be used to establish the
real-rootedness of dual Chow polynomials.

Let $f,g\in \mathbb{R}[x]$ be two real-rooted polynomials with positive leading coefficients. 
Following Br{\"a}nd{\'e}n~\cite{Bra15}, we say that $f$ is an \emph{interleaver} of $g$ (written $f \preceq g$) if
$$
\cdots \le r_2(f) \le r_2(g)\le r_1(f)\le r_1(g),
$$
where $r_i(f)$ and $r_i(g)$ are the zeros of $f$ and $g$ in non-increasing order, respectively. 
For notational convenience, let $0 \preceq 0$, $0 \preceq h$ and $h \preceq 0$, where $h$ is any real-rooted polynomial with positive leading coefficient. 
By definition, if $f\preceq g$, then $\deg(g)=\deg(f)$ or $\deg(f)+1$. 
A sequence $(f_0, f_1, \dots , f_n)$ of real-rooted polynomials is called \emph{interlacing} if $f_i \preceq f_j$ for all $0 \le i < j \le n$.
We repeatedly use the following basic properties in this paper.

\begin{lemma}[{\cite[Section~7]{Bra15}}]
\label{lem:closure}
Let $f_1,\ldots,f_n,g$ be nonzero real-rooted polynomials with nonnegative coefficients, and let $f=\sum_{i=1}^n a_i f_i$, where $a_1,\ldots,a_n$ are nonnegative and not all zero.
Then the following statements hold:
\begin{enumerate}
\item[(i)] If $g\preceq f_i$ for every $i\in[n]$,
then $f$ is real-rooted and $g\preceq f$.
\item[(ii)] If $f_i\preceq g$ for every $i\in[n]$,
then $f$ is real-rooted and $f\preceq g$.
\item[(iii)] For each $i\in[n]$,
$g\preceq f_i$ if and only if $f_i\preceq xg$.
\end{enumerate}
\end{lemma}

\begin{lemma}[{\cite[Theorem 7.8.5]{Bra15}}]
\label{lem:upper-sum}
Let $(f_0,\ldots,f_n)$ be an interlacing sequence of polynomials with nonnegative coefficients.
Then
\[
\sum_{i=0}^n f_i
\preceq
\sum_{i=k}^n f_i
\]
for every $0\leq k\leq n$.
\end{lemma}

The following result relates the interlacing of
$\gamma$-positive polynomials to that of their
associated $\gamma$-polynomials.

\begin{proposition}[{\cite[Lemma~4.10]{HJM19}}]
\label{prop:gamma-zeros}
Let $f(x)$ and $g(x)$ be $\gamma$-positive polynomials
with centers of symmetry $d/2$ and $(d+1)/2$,
respectively, and suppose that $f(0),g(0)>0$.
If
\[
\gamma(f;z)\preceq\gamma(g;z),
\]
then $f(x)\preceq g(x)$.
\end{proposition}

In view of Proposition~\ref{prop:gamma-zeros}, to prove Theorem~\ref{thm:main} for matroids of positive rank, it suffices to establish the corresponding interlacing relation between the $\gamma$-polynomials of $\Hd_M(x)$ and $\Fd_M(x)$ and verify the positivity of their constant terms.
\section{Decomposition formulas for dual Chow polynomials}
\label{sec:decomposition}
In this section, we express the dual Chow and dual right augmented Chow polynomials of paving matroids in terms of those of uniform matroids.
Throughout this section and the next, unless otherwise stated, let $M$ be a paving matroid of rank $r\geq 2$
on an $n$-element ground set $E$.
Since $M$ is paving, every subset of $E$ of size at most
$r-1$ is independent.
It follows that the flats of rank at most $r-2$ are precisely the subsets of $E$ of size at most $r-2$.
In particular, $M$ is loopless.

A hyperplane $H$ of $M$ is called \emph{nontrivial} if $|H|\ge r$.
Let
\[
\cH(M)=\{H\in \mathcal{L}(M):\rk(H)=r-1,\ |H|\geq r\}.
\]
denote the set of nontrivial hyperplanes of $M$.

Since every $(r-1)$-element subset of $E$ is contained in a unique hyperplane, we have
\begin{equation}\label{eq:packing}
\sum_{H\in\cH(M)}\binom{|H|}{r-1}\leq\binom n{r-1}.
\end{equation}

We first compute $\mu_M(F,E)$ for flats $F$ of rank at most $r-2$.

\begin{lemma}\label{lem:mobius}
Let $M$ be a paving matroid of rank $r\ge 2$ on an $n$-element
ground set $E$. If $F$ is a flat of $M$ with $|F|=k\le r-2$, then
\begin{equation}\label{eq:mobius}
\mu_M(F,E)=(-1)^{r-k}\binom{n-k-1}{r-k-1}
+(-1)^{r-1-k}\!\sum_{\substack{H\in\cH(M)\\ F\subseteq H}}
\binom{|H|-k-1}{r-k-1}.
\end{equation}
In particular,
\begin{equation}\label{eq:mobius-constant}
(-1)^r\mu_M(\varnothing,E)=\binom{n-1}{r-1}
-\sum_{H\in\cH(M)}\binom{|H|-1}{r-1}.
\end{equation}
\end{lemma}
\begin{proof}
Let $\U_{r,n}$ be the uniform matroid of rank $r$ on $E$.
For $N\in\{M,\U_{r,n}\}$, the defining recurrence for the
M\"obius function gives
\[
\mu_N(F,E)
=-\sum_{\substack{\tilde{F}\in\mathcal L(N)\\ F\subseteq \tilde{F}\subsetneq E}}
\mu_N(F,\tilde{F}).
\]
Since $M$ is paving, every subset of $E$ of size at most $r-1$ is independent. 
If $|\tilde{F}|\le r-2$, then $\tilde{F}\cup\{e\}$ is independent for every $e\in E\setminus \tilde{F}$. 
Thus adjoining any element outside $\tilde{F}$ increases its rank by one, and $\tilde{F}$ is a flat. 
Conversely, any subset of size at least $r-1$ contains an independent $(r-1)$-element subset and therefore has rank at least $r-1$.
It follows that the flats of rank at most $r-2$ in both $M$ and $\U_{r,n}$ are precisely the subsets of $E$ of size at most $r-2$.
For each such flat $\tilde{F}$ containing $F$, the interval $[F,\tilde{F}]$ is Boolean in both lattices. 
Hence the terms indexed by these flats agree in the two recurrences.

We first compute $\mu_{\U_{r,n}}(F,E)$. The proper flats of
$\U_{r,n}$ are precisely the subsets of $E$ of size at most $r-1$.
For $k\le i\le r-1$, there are $\binom{n-k}{i-k}$ such flats $\tilde{F}$
of size $i$ containing $F$. Since $[F,\tilde{F}]$ is a Boolean lattice of
rank $i-k$, we have $\mu_{\U_{r,n}}(F,\tilde{F})=(-1)^{i-k}$.
Consequently,
\[
\mu_{\U_{r,n}}(F,E)
=-\sum_{i=k}^{r-1}(-1)^{i-k}\binom{n-k}{i-k}
=(-1)^{r-k}\binom{n-k-1}{r-k-1},
\]
where the last equality follows from the identity
\[
\sum_{i=0}^{s}(-1)^i\binom{n}{i}
=(-1)^s\binom{n-1}{s},
\qquad 0\le s<n.
\]

Next, let $I$ be an $(r-1)$-element subset of $E$ containing $F$.
Since $I$ is independent, its closure in $M$ is a hyperplane.
Any hyperplane containing $I$ also contains its closure and has the
same rank, so it equals the closure of $I$. We may therefore partition
the $(r-1)$-element subsets containing $F$ according to their
closures in $M$. If $I$ is itself a hyperplane of $M$, then it is
the only $(r-1)$-element subset with that closure. The interval
$[F,I]$ is Boolean in both lattices, so the term indexed by $I$
is the same in the two recurrences.

Now fix a nontrivial hyperplane $H$ containing $F$. 
The $(r-1)$-element subsets containing $F$ whose closure
is $H$ are exactly those contained in $H$. 
There are
$\binom{|H|-k}{r-1-k}$ of them, and the values of the M\"obius function in
$\U_{r,n}$ sum to $(-1)^{r-1-k}\binom{|H|-k}{r-1-k}$. 
The corresponding summand
in the recurrence for $M$ is $\mu_M(F,H)$.

Every flat properly contained in $H$ has rank at most $r-2$.
Thus the interval $[F,H)$ consists of the sets $\tilde{F}$ with
$F\subseteq \tilde{F}\subseteq H$ and $|\tilde{F}|\le r-2$. For each
$k\le i\le r-2$, there are $\binom{|H|-k}{i-k}$ such flats of
size $i$, and $\mu_M(F,\tilde{F})=(-1)^{i-k}$ for each one. Applying
the M\"obius function on $[F,H]$ gives
\[
\mu_M(F,H)
=-\sum_{i=k}^{r-2}(-1)^{i-k}\binom{|H|-k}{i-k}
=(-1)^{r-1-k}\binom{|H|-k-1}{r-2-k}.
\]

By the recurrence for the M\"obius function, the contribution of $H$ 
to $\mu_M(F,E)-\mu_{\U_{r,n}}(F,E)$ is given by
\[
\sum_{\substack{F\subseteq\tilde{F}\subseteq H\\
|\tilde{F}|=r-1}}
\mu_{\U_{r,n}}(F,\tilde{F})-\mu_M(F,H).
\]
We thus obtain
\begin{align*}
(-1)^{r-1-k}\binom{|H|-k}{r-1-k}-\mu_M(F,H)
&=(-1)^{r-1-k}
\left(
\binom{|H|-k}{r-1-k}-\binom{|H|-k-1}{r-2-k}
\right)\\
&=(-1)^{r-1-k}\binom{|H|-k-1}{r-1-k}.
\end{align*}
Summing over the nontrivial hyperplanes containing $F$ and adding $\mu_{\U_{r,n}}(F,E)$ yields~\eqref{eq:mobius}.
Setting $F=\varnothing$ and multiplying by $(-1)^r$ gives~\eqref{eq:mobius-constant}.
This completes the proof.
\end{proof}

We next recall the definition of Eulerian polynomials; see~\cite{Pet15}.
For $n\geq1$, write $\mathfrak S_n$ for the set of permutations of $[n]$.
The \emph{descent set} of $w\in\mathfrak S_n$ is
\[
\Des(w)=\{i\in[n-1]:w(i)>w(i+1)\}.
\]
Set $\des(w)=|\Des(w)|$.
The \emph{Eulerian polynomials} are given by $A_0(x)=1$ and
\[
A_n(x)=\sum_{w\in\mathfrak S_n}x^{\des(w)},
\quad n\geq1.
\]

For an integer $m\geq1$, set
\begin{equation*}
Q_m(x)=\frac{x^m-x}{x-1}=\sum_{i=1}^{m-1}x^i,
\qquad Q_1(x)=0.
\end{equation*}
Let $P$ be a finite bounded graded poset of rank $r\geq 1$.
By~\cite[Theorem~4.1]{CFH26}, we have
\begin{equation}\label{eq:chain}
\Hd_P(x)=(-1)^r
\sum_{\zmin\leq c_0<c_1<\cdots<c_s=\zmax}
\mu_P(\zmin,c_0)
\prod_{i=1}^s\mu_P(c_{i-1},c_i)Q_{\rk(c_i)-\rk(c_{i-1})}(x).
\end{equation}
The term for $s=0$ is $(-1)^r\mu_P(\zmin,\zmax)$. 
We take the dual Chow polynomial of a one-element interval to be $1$.
Grouping the remaining chains according to $F=c_{s-1}$ and applying~\eqref{eq:chain} to $[\zmin,F]$, we obtain

\begin{equation}\label{eq:last-flat}
\Hd_P(x)=(-1)^r\mu_P(\zmin,\zmax)
+\sum_{F<\zmax}(-1)^{r+\rk(F)}
\Hd_{[\zmin,F]}(x)\mu_P(F,\zmax)Q_{r-\rk(F)}(x).
\end{equation}

\begin{theorem}\label{thm:correction}
If $M$ is a paving matroid of rank $r\geq2$ on $n$ elements, then
\begin{equation}\label{eq:correction}
\Hd_M(x)=\Hd_{\U_{r,n}}(x)
-\sum_{H\in\cH(M)}\Hd_{\U_{r,|H|}}(x).
\end{equation}
\end{theorem}
\begin{proof}
Let $\U_{r,n}$ be the uniform matroid on $E$. 
Since $Q_1(x)=0$, the sum in~\eqref{eq:last-flat} reduces to the sum over flats of rank at most $r-2$.
These flats are the same in $M$ and $\U_{r,n}$.
For each such flat $F$ of rank $k$, the interval $[\varnothing,F]$ is a Boolean lattice of rank $k$.
By~\cite[Proposition~6.8]{CFH26}, we have $\Hd_{[\varnothing,F]}(x)=A_k(x)$.
Hence, for $N\in\{M,\U_{r,n}\}$, we have
\[
\Hd_N(x)
=(-1)^r\mu_N(\varnothing,E)
+\sum_{k=0}^{r-2} A_k(x)Q_{r-k}(x)
\sum_{\substack{I\subseteq E\\ |I|=k}}
(-1)^{r+k}\mu_N(I,E).
\]
By Lemma~\ref{lem:mobius}, for $F\subseteq E$ with $|F|=k\leq r-2$, we have
\[
(-1)^{r+k}\bigl(\mu_{\U_{r,n}}(F,E)-\mu_M(F,E)\bigr)
=\sum_{\substack{H\in\cH(M)\\ F\subseteq H}}
\binom{|H|-k-1}{r-k-1}.
\]

Taking $F=\varnothing$ in the above identity, we see that each $H\in\cH(M)$ contributes $\binom{|H|-1}{r-1}$ to the difference of the constant terms.
Since every $k$-element subset of $H$ is a flat of $M$ for $0\leq k\leq r-2$, there are exactly $\binom{|H|}{k}$ flats $F$ of rank $k$ contained in $H$.
Thus, in the sum over flats of rank $k$, the term $\binom{|H|-k-1}{r-k-1}$ corresponding to $H$ occurs $\binom{|H|}{k}$ times.
Subtracting the expression for $\Hd_M(x)$ from that for $\Hd_{\U_{r,n}}(x)$ and interchanging the order of summation, we obtain
\[
\Hd_{\U_{r,n}}(x)-\Hd_M(x)
=\sum_{H\in\cH(M)}
\left[
\binom{|H|-1}{r-1}
+\sum_{k=0}^{r-2}
\binom{|H|}{k}\binom{|H|-k-1}{r-k-1}
A_k(x)Q_{r-k}(x)
\right].
\]
Each $H\in\cH(M)$ has size at least $r$. By the uniform matroid
formula~\cite[Proposition~6.9]{CFH26}, the expression in brackets is
$\Hd_{\U_{r,|H|}}(x)$.
This proves~\eqref{eq:correction}.
\end{proof}

Let $P$ be a finite bounded graded poset of rank $r\geq 1$.
By~\eqref{eq:intro-flag-H} and~\eqref{eq:intro-flag-F}, the associated $\gamma$-polynomials are given by
\begin{align}
\Gamma_P(z)
&:=\gamma(\Hd_P;z)
=\sum_{\substack{
S\subseteq[r-1]\\
r-1\notin S\\
S\text{ stable}
}}
\beta_P([r-1]\setminus S)z^{|S|}, \label{eq:flag-H}\\
\Gt_P(z)
&:=\gamma(\Fd_P;z)
=\sum_{\substack{
S\subseteq[r-1]\\
S\text{ stable}
}}
\beta_P([r-1]\setminus S)z^{|S|}.
\label{eq:flag-F}
\end{align}

For matroids, we use the matroid as a subscript in place of its lattice of flats. 
The following lemma relates $\Gamma_M(z)$ and $\Gt_M(z)$ for paving matroids.

\begin{lemma}\label{lem:augmented}
For a paving matroid $M$ of rank $r\geq2$ on $n$ elements, we have
\begin{align}
\Gt_M(z)&=\Gamma_M(z)+z\Gamma_{\U_{r-1,n}}(z),\label{eq:augmented-gamma}\\
\Fd_M(x)&=(1+x)\Hd_M(x)+x\Hd_{\U_{r-1,n}}(x).\label{eq:augmented-x}
\end{align}
\end{lemma}

\begin{proof}
By~\eqref{eq:flag-H} and~\eqref{eq:flag-F}, the difference between the two $\gamma$-polynomials consists of the terms indexed by stable subsets $S\subseteq[r-1]$ containing $r-1$. 
Thus
\[
\Gt_M(z)-\Gamma_M(z)
=
\sum_{\substack{
S\subseteq[r-1]\\
S\text{ stable}\\
r-1\in S
}}
\beta_M([r-1]\setminus S)z^{|S|}.
\]

If $r=2$, the only stable subset of $[1]$ containing $1$ is $\{1\}$.
Since $\beta_M(\varnothing)=1$, we obtain
\[
\Gt_M(z)-\Gamma_M(z)=z.
\]
As $\Gamma_{\U_{1,n}}(z)=1$, this proves~\eqref{eq:augmented-gamma} for $r=2$.
Now let $r\geq3$, since $S$ is stable and contains $r-1$, it excludes $r-2$.
Write $T=S\setminus\{r-1\}$. 
Then $T$ is a stable subset of $[r-3]$, and $[r-1]\setminus S=[r-2]\setminus T$, $|S|=|T|+1$.
Hence
\[
\Gt_M(z)-\Gamma_M(z)
=
z\sum_{\substack{
T\subseteq[r-3]\\
T\text{ stable}
}}
\beta_M([r-2]\setminus T)z^{|T|}.
\]

The subsets of $E$ of size at most $r-2$ are flats of both $M$ and $\mathsf U_{r-1,n}$. 
Thus the subposets of $\mathcal{L}(M)$ and $L(\mathsf U_{r-1,n})$ selected by the ranks $[r-2]$ coincide.
The corresponding flag $h$-vector entries therefore satisfy
\[
\beta_M(A)=
\beta_{\mathsf U_{r-1,n}}(A)
\quad\text{for every }A\subseteq[r-2].
\]
Substituting these flag $h$-vector entries gives
\[
\Gt_M(z)-\Gamma_M(z)
=
z\sum_{\substack{
T\subseteq[r-3]\\
T\text{ stable}
}}
\beta_{\mathsf U_{r-1,n}}([r-2]\setminus T)
z^{|T|}=z\Gamma_{\mathsf U_{r-1,n}}(z).
\]
Together with the case $r=2$, this proves~\eqref{eq:augmented-gamma} for $r\geq2$.

Use the $\gamma$-expansion in Section~\ref{sec:pre}, we have
\begin{align*}
\Hd_M(x)
&=(1+x)^{r-1}
\Gamma_M\left(\frac{x}{(1+x)^2}\right),\\
\Fd_M(x)
&=(1+x)^r
\Gt_M\left(\frac{x}{(1+x)^2}\right),\\
\Hd_{\mathsf U_{r-1,n}}(x)
&=(1+x)^{r-2}
\Gamma_{\mathsf U_{r-1,n}}
\left(\frac{x}{(1+x)^2}\right).
\end{align*}
Substituting~\eqref{eq:augmented-gamma} into the formula for the dual right augmented Chow polynomial gives
\begin{align*}
\Fd_M(x)
&=(1+x)^r
\Gamma_M\left(\frac{x}{(1+x)^2}\right)+x(1+x)^{r-2}
\Gamma_{\mathsf U_{r-1,n}}
\left(\frac{x}{(1+x)^2}\right)\\
&=(1+x)\Hd_M(x)
+x\Hd_{\mathsf U_{r-1,n}}(x).
\end{align*}
This proves~\eqref{eq:augmented-x}.
\end{proof}

\begin{example}
Let $M$ be the rank-$3$ matroid on $E=[6]$ whose bases are all $3$-element subsets of $E$ except $\{1,2,3\}$ and $\{4,5,6\}$.
Then $\cH(M)=\bigl\{\{1,2,3\},\{4,5,6\}\bigr\}$.
By Theorem~\ref{thm:correction}, we have
\[
\Hd_M(x)
=
\Hd_{\U_{3,6}}(x)
-2\Hd_{\U_{3,3}}(x).
\]
Since
\[
\Hd_{\U_{3,6}}(x)=10+34x+10x^2,
\qquad
\Hd_{\U_{3,3}}(x)=1+4x+x^2,
\]
we obtain
\[
\Hd_M(x)=8+26x+8x^2.
\]
Since $\Hd_{\U_{2,6}}(x)=5+5x$, Lemma~\ref{lem:augmented} gives
\[
\Fd_M(x)
=(1+x)\Hd_M(x)+x\Hd_{\U_{2,6}}(x)
=8+39x+39x^2+8x^3.
\]
The two $\gamma$-polynomials are therefore
\[
\Gamma_M(z)=8+10z,
\qquad
\Gt_M(z)=8+15z.
\]
\end{example}

Theorem~\ref{thm:correction} and Lemma~\ref{lem:augmented} imply the following corollary.
\begin{corollary}
If $M$ is a paving matroid of rank $r\geq2$, then
\begin{equation*}\label{eq:F-correction}
\Fd_M(x)=\Fd_{\U_{r,n}}(x)
-(1+x)\sum_{H\in\cH(M)}\Hd_{\U_{r,|H|}}(x).
\end{equation*}
\end{corollary}

Finally, we determine the constant terms and degrees
of $\Hd_M(x)$ and $\Fd_M(x)$.

\begin{lemma}\label{lem:positive-constant}
Let $M$ be a paving matroid of rank $r\geq 2$ on $n$ elements.
Then
\[
\Hd_M(0)=\Fd_M(0)
=\binom{n-1}{r-1}
-\sum_{H\in\cH(M)}\binom{|H|-1}{r-1}>0.
\]
Moreover, $\deg\Hd_M=r-1$ and $\deg\Fd_M=r$.
\end{lemma}

\begin{proof}
The formula for $\Hd_M(0)$ follows from~\eqref{eq:correction} and~\cite[Proposition~6.9]{CFH26}.
By Lemma~\ref{lem:augmented}, we have $\Fd_M(0)=\Hd_M(0)$.

If $\cH(M)=\varnothing$, this term is positive. Otherwise, for
each $H\in\cH(M)$, the strict inequality below follows from
$r\geq2$ and $|H|<n$.
\[
\binom{|H|-1}{r-1}
=\frac{|H|-r+1}{|H|}
\binom{|H|}{r-1}
<
\frac{n-r+1}{n}
\binom{|H|}{r-1}.
\]
Summing over $H\in\cH(M)$ and using~\eqref{eq:packing} gives
\begin{align*}
\sum_{H\in\cH(M)}
\binom{|H|-1}{r-1}
&<
\frac{n-r+1}{n}
\sum_{H\in\cH(M)}
\binom{|H|}{r-1}\\
&\leq
\frac{n-r+1}{n}\binom{n}{r-1}
=\binom{n-1}{r-1}.
\end{align*}
The common constant term is therefore positive.

By~\eqref{eq:intro-flag-H} and~\eqref{eq:intro-flag-F}, the centers of symmetry of $\Hd_M(x)$ and $\Fd_M(x)$ are $(r-1)/2$ and $r/2$, respectively.
Since their constant terms are positive, we have $\deg\Hd_M=r-1$ and $\deg\Fd_M=r$.
\end{proof}

\section{The proof of main theorem}
\label{sec:eulerian}
In this section, we prove Theorem~\ref{thm:main} using the decomposition formulas established in the preceding section and the interlacing properties of refined Eulerian polynomials.

Let $n\geq 1$ and $0\leq k\leq n$. Define
\[
\widehat{\mathfrak S}_{n+1,k+1}
=
\{w\in\mathfrak S_{n+1}:w(1)=k+1,\ \Des(w)\text{ is stable}\}.
\]
Following Hoster and Stump~\cite{HS26}, for $S\subseteq T\subseteq[n]$, 
we define the refined Eulerian polynomials by
\[
p_{n,k}^{S\subseteq T}(x)
=
\sum_{\substack{
w\in\widehat{\mathfrak S}_{n+1,k+1}\\
S\subseteq\Des(w)\subseteq T
}}
x^{\des(w)}.
\]
We write $p_{n,k}^T(x)$ for $p_{n,k}^{\varnothing\subseteq T}(x)$, with the conventions $p_{0,0}^{\varnothing}(x)=1$ and $p_{0,0}^{\{1\}\subseteq[0]}(x)=0$.

Since every permutation in $\mathfrak S_{n+1}$ beginning with $n+1$ has a descent at position $1$, whereas one beginning with $1$ does not, we have
$p_{n,n}^{[n]\setminus\{1\}}(x)=p_{n,0}^{\{1\}\subseteq[n]}(x)=0$ for any $n\ge1$.
According to whether $1\in\Des(w)$, we obtain
\[
p_{n,i}^{[n]}(x)=p_{n,i}^{[n]\setminus\{1\}}(x)+p_{n,i}^{\{1\}\subseteq[n]}(x),\quad \text{for } 0\le i\le n.
\]

We will use the following interlacing result of Hoster and Stump.

\begin{lemma}[{\cite[Theorem~3.3]{HS26}}]
\label{lem:HS-interlacing}
Let $n\geq 2$ and $T\in\{[n],[n-1]\}$.
Then the sequences
\[
\left(
p_{n,0}^{T\setminus\{1\}},\ldots,
p_{n,n}^{T\setminus\{1\}},
p_{n,0}^{\{1\}\subseteq T},\ldots,
p_{n,n}^{\{1\}\subseteq T}
\right)
\]
and
\[
\left(p_{n,0}^{T},p_{n,1}^{T},\ldots,p_{n,n}^{T}\right)
\]
are interlacing. Moreover,
\[
p_{n,i}^{T\setminus\{1\}}\preceq p_{n,j}^{T}
\quad\text{and}\quad
p_{n,i}^{T}\preceq p_{n,j}^{\{1\}\subseteq T}
\]
whenever $0\leq i\leq j\leq n$.
\end{lemma}

Then we have the following lemma.

\begin{lemma}
\label{lem:eulerian-input}
For every $n\geq0$, the sequence
\begin{equation}\label{eq:nonzero-path}
\left(p_{n,0}^{[n]\setminus\{1\}},\ldots,p_{n,n}^{[n]\setminus\{1\}},p_{n,0}^{\{1\}\subseteq[n]},\ldots,p_{n,n}^{\{1\}\subseteq[n]}\right)
\end{equation}
is interlacing.
\end{lemma}
\begin{proof}
For $n=0,1$, the sequences in~\eqref{eq:nonzero-path} are $(1,0)$ and $(1,0,0,x)$, respectively, and hence are interlacing.
For $n\geq 2$, the result follows from~Lemma~\ref{lem:HS-interlacing} by taking $T=[n]$ and including the two zero polynomials $p_{n,n}^{[n]\setminus\{1\}}(x)$ and $p_{n,0}^{\{1\}\subseteq[n]}(x)$.
\end{proof}

Caiolo, Ferroni, and Hoster~\cite[Proposition~6.10]{CFH26} obtained the following expansion for the $\gamma$-polynomial of a uniform matroid:
\begin{equation}\label{eq:uniform-gamma}
\Gamma_{\U_{r,n}}(z)=\sum_{k=0}^{r-1}
\binom{n-1-k}{r-1-k}p_{r-1,k}^{[r-1]\setminus\{1\}}(z).
\end{equation}
By deleting the first entry of a permutation, we obtain the following identity.

\begin{lemma}\label{lem:deletion}
For $r\geq2$ and $0\leq k\leq r-1$, we have
\begin{equation}\label{eq:deletion}
p_{r-1,k}^{[r-1]\setminus\{1\}}(x)
=\sum_{i=k}^{r-2}p_{r-2,i}^{[r-2]}(x).
\end{equation}
\end{lemma}

\begin{proof}
For $k=r-1$, both sides of~\eqref{eq:deletion} are zero.
Suppose that $0\leq k\leq r-2$.

Let $w$ be a permutation counted by $p_{r-1,k}^{[r-1]\setminus\{1\}}(x)$.
Then $w(1)=k+1$ and $w(2)>k+1$.
Deleting the first entry and replacing the remaining entries by $1,\ldots,r-1$ in increasing order, we obtain a permutation $v\in\mathfrak S_{r-1}$.
Since $v(1)=w(2)-1$, we have $v(1)=i+1$ for some $k\leq i\leq r-2$.
Moreover,
\[
\Des(v)=\{j-1:j\in\Des(w)\}.
\]
Thus $\Des(v)$ is stable and $\des(v)=\des(w)$.
It follows that $v$ is counted by $p_{r-2,i}^{[r-2]}(x)$.

Conversely, for a permutation counted by $p_{r-2,i}^{[r-2]}(x)$ with $i\ge k$, increase each entry greater than $k$ by one.
This gives a permutation of $[r]\setminus\{k+1\}$ whose first entry is $i+2>k+1$. 
Insert $k+1$ at the beginning. 
This creates no descent at position $1$ and shifts each of the other descents one position to the right. 
The resulting permutation is counted by $p_{r-1,k}^{[r-1]\setminus\{1\}}(x)$.
These two operations are inverse to each other and
preserve the number of descents.
This completes the proof.
\end{proof}

Combining Theorem~\ref{thm:correction} with~\eqref{eq:uniform-gamma} and~\eqref{eq:deletion}, we obtain the following expansion for $\Gamma_M(z)$.

\begin{proposition}\label{prop:common}
Let $M$ be a paving matroid of rank $r\geq2$ on $n$ elements. Then
\begin{equation}\label{eq:common}
\Gamma_M(z)
=\sum_{i=0}^{r-2}
\left[
\sum_{k=0}^{i}
\left(
\binom{n-1-k}{r-1-k}
-\sum_{H\in\cH(M)}
\binom{|H|-1-k}{r-1-k}
\right)
\right]\cdot
p_{r-2,i}^{[r-2]}(z).
\end{equation}
\end{proposition}

\begin{proof}
Since all the polynomials in~\eqref{eq:correction} are symmetric about $(r-1)/2$, taking their $\gamma$-polynomials gives
\[
\Gamma_M(z)
=
\Gamma_{\U_{r,n}}(z)
-
\sum_{H\in\cH(M)}\Gamma_{\U_{r,|H|}}(z).
\]
Substituting~\eqref{eq:uniform-gamma} into the above identity and using Lemma~\ref{lem:deletion}, we obtain
\begin{align*}
\Gamma_M(z)
&=
\sum_{k=0}^{r-1}
\left[
\binom{n-1-k}{r-1-k}
-
\sum_{H\in\cH(M)}
\binom{|H|-1-k}{r-1-k}
\right]
p_{r-1,k}^{[r-1]\setminus\{1\}}(z)\\
&=
\sum_{k=0}^{r-1}
\left[
\binom{n-1-k}{r-1-k}
-
\sum_{H\in\cH(M)}
\binom{|H|-1-k}{r-1-k}
\right]
\sum_{i=k}^{r-2}p_{r-2,i}^{[r-2]}(z).
\end{align*}
Interchanging the order of summation yields~\eqref{eq:common}.
\end{proof}

We next compare the weights in the expansions of $\Gamma_M(z)$ and $\Gamma_{\U_{r-1,n}}(z)$.

\begin{proposition}\label{prop:ratios}
Let $M$ be a paving matroid of rank $r\geq2$ on $n$ elements. 
For $0\leq i\leq r-2$, let
\[
u_i=\binom{n-1-i}{r-2-i},\quad \text{and}\quad 
v_i=\sum_{k=0}^{i}\left[
\binom{n-1-k}{r-1-k}
-\sum_{H\in\cH(M)}\binom{|H|-1-k}{r-1-k}\right].
\]
Then
\begin{equation*}
0<\frac{v_0}{u_0}\leq\frac{v_1}{u_1}\leq\cdots\leq\frac{v_{r-2}}{u_{r-2}}.
\end{equation*}
In particular, every $v_i$ is positive.
\end{proposition}

\begin{proof}
By the hockey-stick identity (see~\cite[Eq.~(5.10)]{GKP94}), we have
\[
\sum_{k=0}^i\binom{n-1-k}{r-1-k}
=
\binom{n}{r-1}-\binom{n-1-i}{r-2-i}.
\]
It follows that
\[
v_i
=
\binom{n}{r-1}-u_i
-
\sum_{H\in\cH(M)}
\left[
\binom{|H|}{r-1}
-
\binom{|H|-1-i}{r-2-i}
\right].
\]
Hence
\[
\frac{v_i}{u_i}
=
\frac{\displaystyle
\binom{n}{r-1}
-
\sum_{H\in\cH(M)}\binom{|H|}{r-1}}
{u_i}
-1+
\sum_{H\in\cH(M)}
\frac{\binom{|H|-1-i}{r-2-i}}{u_i}.
\]
By~\eqref{eq:packing}, we obtain that
$\binom{n}{r-1}-\sum_{H\in\cH(M)}\binom{|H|}{r-1}\ge 0$.
Moreover, for $0\leq i<r-2$, we have
\[
\frac{u_{i+1}}{u_i}
=
\frac{r-2-i}{n-1-i}<1.
\]
Thus $1/u_i$ is increasing in $i$.
Since $r\leq |H|<n$ for every $H\in\cH(M)$, we also have
\begin{align*}
\frac{\binom{|H|-2-i}{r-3-i}}{u_{i+1}}
&=
\frac{n-1-i}{|H|-1-i}
\frac{\binom{|H|-1-i}{r-2-i}}{u_i}
>
\frac{\binom{|H|-1-i}{r-2-i}}{u_i}.
\end{align*}
Hence $v_i/u_i$ is nondecreasing in $i$, that is,
\[
\frac{v_0}{u_0}
\leq
\frac{v_1}{u_1}
\leq\cdots\leq
\frac{v_{r-2}}{u_{r-2}}.
\]

By Lemma~\ref{lem:positive-constant}, we have
\[
v_0
=
\binom{n-1}{r-1}
-
\sum_{H\in\cH(M)}
\binom{|H|-1}{r-1}
=
\Hd_M(0)>0.
\]
Since $u_i>0$ for all $0\leq i\leq r-2$, the desired inequalities follow.
\end{proof}

To compare the above two expansions, we establish the following lemma on interlacing polynomials.

\begin{lemma}\label{lem:weighted}
Let $g_0,\ldots,g_n,f_0,\ldots,f_n$ be polynomials with nonnegative
coefficients satisfying the following conditions.
\begin{enumerate}
\item[(i)] $(g_0,\ldots,g_n)$ is an interlacing sequence, and $g_i\not\equiv 0$ for at least one $i$.
\item[(ii)] $g_i\preceq f_j$ for all $0\le i,j\le n$.
\end{enumerate}
Let $a_i\ge 0$ and $b_i>0$ for $0\le i\le n$. If
\[
\frac{a_0}{b_0}\le \frac{a_1}{b_1}\le \cdots \le \frac{a_n}{b_n},
\]
then
\[
\sum_{i=0}^n b_ig_i
\preceq
\sum_{i=0}^n a_i(f_i+g_i).
\]
\end{lemma}

\begin{proof}
Let $h_i=f_i+g_i$. 
In applying Lemma~\ref{lem:closure}, we omit terms with zero coefficient or zero polynomial. 
Since $(g_0,\ldots,g_n)$ is an interlacing sequence and $b_i>0$ for every $i$, the sequence $(b_0g_0,\ldots,b_ng_n)$ is also interlacing.
By Lemma~\ref{lem:upper-sum}, we have
\begin{equation}\label{eq:partial-g}
\sum_{i=0}^n b_i g_i \preceq\sum_{i=k}^n b_ig_i,
\quad \text{for}\ 0\le k\le n.
\end{equation}

By assumption~(ii) and Lemma~\ref{lem:closure} (ii),
we have $\sum_{i=0}^n b_i g_i \preceq f_j$ for every $0\leq j\leq n$.
Hence Lemma~\ref{lem:closure} (i) gives
\begin{equation}\label{eq:partial-f}
\sum_{i=0}^n b_i g_i \preceq\sum_{i=k}^n b_if_i,
\quad \text{for}\ 0\le k\le n.
\end{equation}
Combining \eqref{eq:partial-g} and \eqref{eq:partial-f},
we obtain
\[
\sum_{i=0}^n b_i g_i \preceq\sum_{i=k}^n b_ih_i,
\quad \text{for}\ 0\le k\le n.
\]

Finally, note that
\[
\sum_{i=0}^n a_ih_i
=
\frac{a_0}{b_0}\sum_{i=0}^n b_ih_i
+
\sum_{k=1}^n
\left(
\frac{a_k}{b_k}-\frac{a_{k-1}}{b_{k-1}}
\right)
\sum_{i=k}^n b_ih_i.
\]
Since
\[
0\leq\frac{a_0}{b_0}
\leq\frac{a_1}{b_1}
\leq\cdots\leq\frac{a_n}{b_n},
\]
the right-hand side is a nonnegative linear combination of polynomials interlaced by $\sum_{i=0}^n b_i g_i$.
By Lemma~\ref{lem:closure} (i), we conclude that
\[
\sum_{i=0}^n b_ig_i
\preceq
\sum_{i=0}^n a_i(f_i+g_i).
\]
This completes the proof.
\end{proof}

Combining Propositions~\ref{prop:common} and~\ref{prop:ratios} with Lemma~\ref{lem:weighted}, we obtain the following theorem.

\begin{theorem}\label{thm:gamma-interlacing}
Let $M$ be a paving matroid of rank $r\geq2$ on $n$ elements. Then
\begin{equation}\label{eq:gamma-interlacing}
\Gamma_{\U_{r-1,n}}(z)\preceq\Gamma_M(z)
\preceq z\Gamma_{\U_{r-1,n}}(z).
\end{equation}
All three polynomials are nonzero and real-rooted.
\end{theorem}

\begin{proof}
If $r=2$, then $\Gamma_{\U_{1,n}}(z)=1$.
Since $\Hd_M(x)$ is symmetric about $1/2$, its $\gamma$-polynomial $\Gamma_M(z)$ is constant and positive by Lemma~\ref{lem:positive-constant}.
Hence~\eqref{eq:gamma-interlacing} holds.

Let $u_i$ and $v_i$ be as in Proposition~\ref{prop:ratios}.
By~\eqref{eq:uniform-gamma} and~\eqref{eq:common}, we have
\[
\Gamma_{\U_{r-1,n}}(z)=\sum_{i=0}^{r-2}u_i\,p_{r-2,i}^{[r-2]\setminus\{1\}}(z),
\]
and 
\[
\Gamma_M(z)=\sum_{i=0}^{r-2}v_i\,p_{r-2,i}^{[r-2]}(z)
=\sum_{i=0}^{r-2}v_i\left(p_{r-2,i}^{[r-2]\setminus\{1\}}(z)+p_{r-2,i}^{\{1\}\subseteq[r-2]}(z)\right).
\]
Set
\[
g_i(z)=p_{r-2,i}^{[r-2]\setminus\{1\}}(z)
\quad\text{and}\quad
f_i(z)=p_{r-2,i}^{\{1\}\subseteq[r-2]}(z),
\quad \text{for}\ 0\leq i\leq r-2.
\]
By Lemma~\ref{lem:eulerian-input}, $(g_0,\ldots,g_{r-2})$ is an interlacing sequence and $g_i\preceq f_j$ for all $0\leq i,j\leq r-2$.
Moreover, $g_0(0)=1$.
By Proposition~\ref{prop:ratios},
\[
0<
\frac{v_0}{u_0}
\le\frac{v_1}{u_1}
\le\cdots\leq\frac{v_{r-2}}{u_{r-2}}.
\]
Applying Lemma~\ref{lem:weighted} with $a_i=v_i$ and $b_i=u_i$, we obtain
\[
\Gamma_{\U_{r-1,n}}(z)\preceq\Gamma_M(z).
\]
Since $\Gamma_{\U_{r-1,n}}(0)=u_0>0$ and $\Gamma_M(0)=\Hd_M(0)>0$ by Lemma~\ref{lem:positive-constant}, both polynomials are nonzero.
The second interlacing relation
\[
\Gamma_M(z)\preceq z\Gamma_{\U_{r-1,n}}(z)
\]
follows from Lemma~\ref{lem:closure} (iii).
This completes the proof.
\end{proof}

We are now prove Theorem~\ref{thm:main}.

\begin{proof}[\textbf{Proof of Theorem~\ref{thm:main}}]
Let $M$ be a finite paving matroid on a ground set $E$ with $|E|=n$.
Write $r=\rk(M)$. 
First assume that $r\ge 2$.
By Theorem~\ref{thm:gamma-interlacing}, we have
\[
\Gamma_M(z)\preceq z\Gamma_{\U_{r-1,n}}(z).
\]
Recall from~\eqref{eq:augmented-gamma} that
\[
\Gt_M(z)=\Gamma_M(z)+z\Gamma_{\U_{r-1,n}}(z).
\]
Since $\Gamma_M(z)$ interlaces both summands on the
right-hand side, Lemma~\ref{lem:closure} (i) gives
\[
\Gamma_M(z)\preceq\Gt_M(z).
\]
Recall that $\Hd_M(x)$ and $\Fd_M(x)$ are $\gamma$-positive with centers of symmetry $(r-1)/2$ and $r/2$, respectively.
Moreover, Lemma~\ref{lem:positive-constant} gives $\Hd_M(0)=\Fd_M(0)>0$.
It follows from Proposition~\ref{prop:gamma-zeros} that
\[
\Hd_M(x)\preceq\Fd_M(x).
\]
Since both polynomials have nonnegative coefficients
and positive constant terms, all their zeros are negative.

If $r=1$, then $\Hd_M(x)=1$ and $\Fd_M(x)=1+x$ by~\eqref{eq:intro-flag-H} and~\eqref{eq:intro-flag-F}.
If $r=0$, then $\Hd_M(x)=\Fd_M(x)=1$.
The result clearly holds in both cases.
This completes the proof.
\end{proof}

\section*{Declaration of AI usage}
During the preparation of this manuscript, the authors used generative AI tools to assist in exploring possible approaches, checking technical details, and improving the exposition. 
The authors verified the mathematical arguments, computations, proofs, and results and take full responsibility for the content of this manuscript.


\end{document}